Stability and bifurcation properties of the algorithms for keeping of differential equations solutions on the required level.

Yu.V. Troshchiev [a]

[a] Dept. Calculate Mathematics and Cybernetics, M.V. Lomonosov Moscow State University, Moscow, 119991, Russia.

Yuri V. Troshchiev, Ph.D., Senior Researcher.

e-mail: yuvt@yandex.ru.

**Abstract.**

Algorithms of control of differential equations solutions are under investigation in the article. Idealized and real modifications of the algorithms are distinguished. An equation, which can be the base equation for investigation of the idealized algorithms properties, is constructed. The difference appearing for real systems and real algorithms is for separate investigation. This difference tends to zero under tending to zero of the time step of control. If the systems of equations satisfy or almost satisfy some properties for which the algorithms are intended, then the results are similar numerically as well. One of the algorithms demonstrates high reliability. Another one is of more complex properties. Bifurcations, periodic solutions and strange attractors are possible in both algorithms in addition to stable steady states.

**Introduction.**

Algorithmically defined methods are often worth applying for keeping different systems states on a necessary level. It is so, for example, for linear systems without steady states in the vicinity of the required state. Thus, appropriate algorithms and investigation of their properties are necessary. Two algorithms (and some modifications) were proposed in [1,2] for control of the atomic reactor power [1-4]. This article investigates properties of these algorithms from stability and the bifurcation properties point of view. Mathematically, it means investigation of a fixed point of some map.

The algorithms are constructed here in the form all the same for what objects they are applied. Thus, they are probably applicable for controlling of different objects under an appropriate choice of parameters, but not only for the fast reactors in the self-adjusting mode.

**1. An object description from the point of view of the algorithms of control.**

Let us describe the system in an idealized form firstly. Let the vector $\mathbf{x} = (x_1, x_2, ..., x_n)$ and time $t$ determine the object state, and the system of differential equations

$$\frac{d\mathbf{x}}{dt} = \mathbf{f}(\mathbf{x}, t). \tag{1}$$

determines the dependence of $\mathbf{x}$ on time. More generally, some components of $\mathbf{x}$ may be numbers and some may be functions. Then the system of equations contains the partial derivatives. Let the first component of $\mathbf{x}$ is some quantity $Q > 0$, which we need to keep approximately equal to $\overline{Q}$. If $Q$ is a functional on system state and is not in a set of variables, then, applying time derivatives of the variables, we can write for $Q$ an equation $Q' = f_1(...)$ and add it to the system. Let, just the same, the second component of $\mathbf{x}$ is some quantity $N$, which is used to control the object. If $N$ is an independent parameter, then we can formally write an equation $N' = 0$ and add it to the system of equations. I.e. we consider that the system (1) is a system with already added, if necessary, the described auxiliary equations. We need it only for simplification of denotation in theoretical investigation not to consider several essentially equal variants. Algorithms are usable without this procedure.

Control over the object is changing of the number $N$ by adding some number $\Delta N$. The control affects not $Q$ directly, but $Q'$, where $Q'$ is the time derivative. We assume that $\lambda = Q'/Q$. The value $\lambda$ is not a component of the vector $\mathbf{x}$. Let the control $\Delta N$ is fulfilled at some moment $t$, then there exists such a function $f$ that

$$\lambda(t + \Delta t) - \lambda(t - 0) = f(\Delta t, \Delta N, \mathbf{x}, t), \ \Delta t \geq 0. \tag{2}$$

The function $f$, generally speaking, depends on the system state $\mathbf{x}$ at the time moment $t$. And $f(0, 0, \mathbf{x}, t) = 0$ always. We also assume that $f$ depends on $\Delta N$ by all means, otherwise, there is no sense to speak about possibility of control.

<u>Proposition.</u> If $f$ does not depends on $\mathbf{x}$ and $t$ then it is linear.

<u>Proof.</u> Really, in this case

$$f(\Delta N, t_2) - f(\Delta N, t_1) = f(0, t_2 - t_1).$$

Dividing by $t_2 - t_1$ and letting $t_2$ go to $t_1$, we obtain

$$f'_t(\Delta N, t_1) = f'_t(0, 0),$$

i.e. the function $f$ can be written in the form $f = b(\Delta N) + \delta_t \Delta t$. Then, linearity of $f$ on $\Delta N$ follows from the equality

$$f(\Delta N_1, 0) + f(\Delta N_2, 0) = f(\Delta N_1 + \Delta N_2, 0).$$

Proof is over.

Let now

$$f(\Delta t, \Delta N) = \delta_t \Delta t + \delta_N \Delta N. \tag{3}$$

Then we can construct a differential equation for which equalities (2), (3) are true:

$$Q' = (\mu + \delta_N N + \delta_t t) Q. \tag{4}$$

We suppose that derivative of the function $Q$ has discontinuity at the moment of the control if $\Delta N \neq 0$: $Q'(t-0) \neq Q'(t+0)$. The system state is defined by three variables: $t$, $Q$ and $N$ here. The equality (3) is true for all values of this variables.

<u>Proposition.</u> If the function $Q(t)$ satisfies the equation (4), then $f(\Delta t, \Delta N)$ is defined by the equality (3). If the function $f(\Delta t, \Delta N)$ is defined by equality (3), then there exists such $\mu$, that $Q(t)$ satisfies the equation (4).

<u>Proof.</u> An expression in parentheses in (4) is $\lambda(t)$. The proof of the first part of the proposition is direct constructing of the function $f(\Delta t, \Delta N)$ for the equation (4):

$$f(\Delta t, \Delta N) = \lambda(t + \Delta t) - \lambda(t) =$$
$$= \mu + \delta_N(N + \Delta N) + \delta_t \Delta t - \mu - \delta_N N = \delta_N \Delta N + \delta_t \Delta t.$$

Let now we know the function $f(\Delta t, \Delta N) = \delta_t \Delta t + \delta_N \Delta N$ and $\lambda(t_0) = \lambda_0$, then $\lambda(t_0 + \Delta t) = \lambda_0 + \delta_N \Delta N + \delta_t \Delta t$. Then, it follows

$$\frac{Q'}{Q} = (\lambda_0 + \delta_N N + \delta_t t) .$$

Proof is over.

The phrase "$Q(t)$ satisfies the equation (4)" does not mean that the system (1) is this equation. It means that the value of $Q(t)$ follows the rule (4). For example, the system

$$\begin{cases} Q' = (\mu + \delta_N N)Q \\ N' = \delta_t / \delta_N \end{cases} \tag{5}$$

forces $Q(t)$ to follow the same rule.

The variable $N$ is an independent parameter in a case of equation (4), thus, it is always possible to choose such an initial value of $N$, that $\mu$ equals 0:

$$Q' = (\delta_N N + \delta_t t)Q . \tag{6}$$

We need this equation in what follows, so, we write its solution here:

$$Q(t) = Q(t_0) \frac{\exp(\delta_N N t + \delta_t t^2 / 2)}{\exp(\delta_N N t_0 + \delta_t t_0^2 / 2)} . \tag{7}$$

The description of the system in this item is called idealized because we suppose possibility of infinitely small values: measurement of instantaneous values of derivatives and instantaneous change of the number $N$. In reality even measurement of $Q$ may be averaging during some time (but we assume that measurement of $Q$ is instantaneous here). The derivative $Q'(t)$ may be measured approximately, for example, by formula $(Q(t) - Q(t-\tau))/\tau$. And control may be done in a time interval. Just the same, there exists some function $f$ that determines change of $\lambda$ similarly to (2). We suppose that there are really measured values in (2) and control takes some time interval in this case (for example, $N$ can change uniformly during the same time interval $\Delta t$). We shall describe it more formally in what follows.

Later we use denotation of subscripts relating to the time $t_i$, like $Q(t_i) = Q_i$, $\lambda(t_i) = \lambda_i$ etc. by default.

**2. Algorithms.**

<u>Algorithm 1.</u> Let $\{t_i\}$, $i = 1,2,...$ is some sequence of time moments $t_{i+1} > t_i$. Actions on odd and even steps are different. We suppose that $\lambda = Q'/Q$, and the control is instantaneous. We suppose also that there is initially given some number $\tilde{\delta}_N$, coinciding with $\delta_N$ by sign and greater than $\delta_N$ in absolute magnitude.

If $i$ is odd, then $N$ is unchanged, i.e. $\Delta N_i = 0$. The time $\Delta t_i = t_{i+1} - t_i$ passed, the value $\tilde{\delta}_t$

$$\frac{\partial f_i}{\partial \Delta t}(0,0) \approx \tilde{\delta}_t = (\lambda_{i+1} - \lambda_i)/\Delta t_i, \tag{8}$$

is calculated, i.e. the formulas (2), (3) at $\Delta N = 0$ are applied, supposing that the linear part of $f$ is enough to be taken into account.

If $i$ is even the desired value of $\lambda$ is calculated by the formula

$$\lambda = \tilde{\lambda}(1 - Q_i/\overline{Q}) \tag{9}$$

and the corresponding change of $N$ by the formula

$$\Delta N_i = (\lambda - \lambda_i - \tilde{\delta}_t \Delta t_i)/\tilde{\delta}_N, \tag{10}$$

i.e. the formulas (2), (3) are secondly applied, supposing that the function $f_i$ differs quite a little from $f_{i-1}$ (and from $f_{i-2}$ at $i \geq 4$). Here $\tilde{\lambda} > 0$ is some parameter of the algorithm, $\tilde{\delta}_t$ is calculated at the previous step by the formula (8). The number $\tilde{\delta}_N$ was given at $i = 2$, and it will be calculated at $i \geq 4$. The number $N$ is to be changed by adding $\Delta N_i$ now. The time $\Delta t_i$ passed, the new value of $\tilde{\delta}_N$ is calculated:

$$\frac{\partial f_i}{\partial \Delta N}(0,0) \approx \tilde{\delta}_N = ((\lambda_{i+1} - \lambda_i) - \tilde{\delta}_t \Delta t_i)/\Delta N, \tag{11}$$

i.e. the formulas (2), (3) are applied the third time.

<u>Modified algorithm 1</u> coincides the algorithm 1, but $\lambda$ now is

$$\lambda_i = \frac{2(Q_i - Q_{i-1})}{\Delta t_{i-1}(Q_i + Q_{i-1})}. \tag{12}$$

$Q_0$ is necessary at the first step in this case. Thus, the time moment $t_0$ is added, and there is no control at the zero step on contrary from the other even steps. As for the other even steps, the value $N$ changes uniformly during $\Delta t_i$ to add $\Delta N_i$ in sum.

<u>Algorithm 2.</u> Let $\{t_i\}$, $i = 1,2,...$ again is some sequence of time moments $t_{i+1} > t_i$. We suppose that $\lambda = Q'/Q$, and control is instantaneous. Let $\tilde{\delta}_N$ is some number such that $\tilde{\delta}_N \delta_N > 0$.

The control is done at each step. The desired $\lambda$ is calculated by the formula (9), then the necessary change of $N$ is calculated by the formula

$$\Delta N_i = \frac{\lambda - \lambda_i}{\tilde{\delta}_N}. \qquad (13)$$

Then $N$ is changed by adding $\Delta N_i$. In contrary from the algorithm 1 $\tilde{\delta}_N$ is one and the same from step to step, and $\Delta N_i$ is calculated neglecting the dependence of the function $f$ on time.

<u>Modified algorithm 2</u> coincides the algorithm 2, but $\lambda$ now is (12) and $N$ changes uniformly during $\Delta t_i$.

### 3. An object for bifurcation analysis.

There are sequences of pairs $\{(Q_i, \lambda_i)\}$ in both algorithms and their modifications, and these pairs follow each other by the relations

$$\begin{pmatrix} Q_{i+1} \\ \lambda_{i+1} \end{pmatrix} = \mathbf{F}_i(Q_i, \lambda_i).$$

If we suppose that the functions $\mathbf{F}_i$ are independent on $i$: $\mathbf{F}_i = \mathbf{F}$, then the problem of investigation of the fixed point and its bifurcations [5-10] for the map $\mathbf{F}$ appears:

$$\begin{pmatrix} Q \\ \lambda \end{pmatrix} = \mathbf{F}(Q, \lambda) = \begin{pmatrix} F_Q(Q, \lambda) \\ F_\lambda(Q, \lambda) \end{pmatrix}. \qquad (14)$$

Independence $\mathbf{F}_i$ on $i$ means independence of composition of maps $\mathbf{F}_{i+1}\mathbf{F}_i$ in the case of the algorithm 1:

$$\begin{pmatrix} Q_{i+2} \\ \lambda_{i+2} \end{pmatrix} = \mathbf{F}_{i+1}(\mathbf{F}_i(Q_i, \lambda_i)). \qquad (15)$$

Taking into account that $\mathbf{F}_i$ corresponds to control of a solution of the system of differential equations, the dependence of $\mathbf{F}_i$ on $i$ means dependence of $\mathbf{F}$ on the other variables of the system (1): $x_2, x_3, \ldots, t$.

<u>Proposition</u>. The function $\mathbf{F}$ does not depends on the system state, excepting $Q$, if and only if the function $f$ is linear and can be written in the form (3).

<u>Proof</u>. Let us prove the proposition for the algorithm 1, for example. Let time $t$ corresponds to an even step. The following relation for the function $Q(t)$ is true

$$(\ln Q(t + 2\Delta t))' = \lambda(t) + f(\Delta N(\Delta t, Q, \lambda, \boldsymbol{\delta}(\mathbf{x}(t - 2\Delta t), t - 2\Delta t)), 2\Delta t, \mathbf{x}, t) =$$
$$= F_\lambda(Q(t), \lambda(t), 2\Delta t, \tilde{\mathbf{x}}, \mathbf{x}(t - 2\Delta t), t).$$

Here $\boldsymbol{\delta}$ is the pair $(\tilde{\delta}_t, \tilde{\delta}_N)$, calculated to the time $t$, $\mathbf{x}(t - 2\Delta t)$ and $t - 2\Delta t$ is the system state at the previous even step.

If $f$ does not depend on $\mathbf{x}$, $t$ and is linear, i.e. is given by (3), then $\boldsymbol{\delta}$ does not depend on $\mathbf{x}(t-2\Delta t)$ and $t-2\Delta t$, thus, $F_\lambda$ does not depend on $\tilde{\mathbf{x}}$, $\mathbf{x}(t-2\Delta t)$, $t$ and $t-2\Delta t$. Then $F_Q$ also does not depends on these values. The second part of the proposition is proved.

It is remained to prove that if $F_\lambda$ does not depends on $\tilde{\mathbf{x}}$, $\mathbf{x}(t-2\Delta t)$ and $t$, then $f$ can be written in the form (3). Let $\boldsymbol{\delta}$ depends on $\mathbf{x}(t-2\Delta t)$ or $t-2\Delta t$, then $F_\lambda$ also depends on $\mathbf{x}(t-2\Delta t)$ or $t-2\Delta t$. The nontrivial dependence of $f$ on $\Delta N$ is recalled here. Thus, $\boldsymbol{\delta}$ does not depends on $\mathbf{x}(t-2\Delta t)$ and $t-2\Delta t$. Let now $f$ nontrivially depends on some component $x_i$, $i \geq 2$. Then there exists such $\Delta t$ and also $\mathbf{x}_1$ and $\mathbf{x}_2$, differing only by values of $x_i$, that

$$f(\Delta N(\Delta t, Q, \lambda, \boldsymbol{\delta}), 2\Delta t, \mathbf{x}_1, t) \neq f(\Delta N(\Delta t, Q, \lambda, \boldsymbol{\delta}), 2\Delta t, \mathbf{x}_2, t).$$

Comparing this inequality with the previous formula for $F_\lambda$ we get a contradiction: the function $F_\lambda$ nontrivially depends on $x_i$. The same is for $t$. So, $f$ can depend only on $Q$.

Now we prove that the dependence $f$ is linear, i.e. it can be written in the form (3). Let us decompose one control action into two ones in two ways:

$$\begin{aligned} f(\Delta N_1 + \Delta N_2, \Delta t, Q_0) &= f(\Delta N_1, 0, Q_0) + f(\Delta N_2, \Delta t, Q_0), \\ f(\Delta N, \Delta t_1 + \Delta t_2, Q_0) &= f(\Delta N, \Delta t_1, Q_0) + f(0, \Delta t_2, Q_1). \end{aligned} \tag{16}$$

In the first equality (16) we subtract the last summand from the both parts of the equation, divide the both parts by $\Delta N_1$ and let $\Delta N_1$ go to zero. We obtain

$$f'_{\Delta N}(\Delta N_2, \Delta t, Q_0) = f'_{\Delta N}(0, 0, Q_0).$$

In the second equality (16) we subtract the first summand of the right part from the both parts of the equation, divide the both parts by $\Delta t_2$ and let $\Delta t_2$ go to zero. We obtain

$$f'_{\Delta t}(\Delta N, \Delta t_1, Q_0) = f'_{\Delta t}(0, 0, Q_1).$$

Now we subtract $f'_{\Delta t}(0,0,Q_0)$ from the both parts, divide them by $\Delta t_1$ and let $\Delta t_1$ go to zero. We obtain $f''_{\Delta t \Delta t}(0,0,Q_0) = f''_{\Delta t Q}(0,0,Q_0)Q'(0) = f''_{\Delta t Q}(0,0,Q_0)\lambda_0 Q_0$. This equality must be true for any $\lambda_0$, including $\lambda_0 = 0$, so,

$$f''_{\Delta t \Delta t}(0,0,Q_0) = 0.$$

But it must be true for $\lambda_0 \neq 0$ too, so,

$$f''_{\Delta t Q}(0,0,Q_0) = 0.$$

It follows from the last four equalities that $f$ can be written in a form $f(\Delta t, \Delta N, Q) = \delta_t \Delta t + \delta_N(Q)\Delta N + \Delta t^3 \psi(\Delta t, Q)$. Comparing the expression for $f(\Delta t, \Delta N, Q)$ with the formulas for calculation of $\boldsymbol{\delta}$ (8), (9), (10), (11), we obtain that $\boldsymbol{\delta}$ does not depend on $\mathbf{x}(t-2\Delta t)$ and $t-2\Delta t$ only if $\delta_N(Q)$ does not depend on $Q$ and $\psi(\Delta t, Q) \equiv 0$. So, the function $f$ can be written in the form (3). Proof is over.

This proposition and the propositions from the item 2 give a possibility to investigate the fixed points and their bifurcations for idealized algorithms applied to the equation (6) and separately investigate what difference appear for applying of real algorithms to real systems of equations.

### 4. Bifurcation properties of the algorithm 1 on the equation (8) with a constant step.

The numbers $\tilde{\delta}_N$ and $\tilde{\delta}_t$ are equal to $\delta_N$ and $\delta_t$ in this case. The cycle of the algorithm consists of two steps, thus, the map (14) is to be constructed on the interval $2\Delta t$. It is convenient to begin from an even step, because the system moves itself during $2\Delta t$ after the control is done. According to (7) and (15) the following formulas define the map (14):

$$Q_{i+2} = Q_i \frac{\exp(\delta_N (N_i + \Delta N)(t_i + 2\Delta t) + \delta_t (t_i + 2\Delta t)^2 / 2)}{\exp(\delta_N (N_i + \Delta N)t_i + \delta_t t_i^2 / 2)} =$$
$$= Q_i \exp\left(2\Delta t \tilde{\lambda}\left(1 - Q_i / \tilde{Q}\right)\right), \qquad (17)$$
$$\lambda_{i+2} = \lambda_i + \delta_N \Delta N + \delta_t 2\Delta t = \tilde{\lambda}\left(1 - Q_i / \tilde{Q}\right) + \delta_t \Delta t,$$

where $\Delta N$ is from the formulas (8), (9), (10), (11). The fixed point of this map exists and unique:

$$\overline{Q} = \tilde{Q},$$
$$\overline{\lambda} = \delta_t \Delta t. \qquad (18)$$

Proposition. If

$$\Delta t < 1/\tilde{\lambda}, \qquad (19)$$

then the fixed point (18) is stable.

Proof. The eigenvalues $\rho_1$, $\rho_2$ of the Jacoby matrix of the map (17) in the point (18) are real:

$$\rho_1 = 1 - 2\Delta t \tilde{\lambda}, \quad \rho_2 = 0.$$

If they are less than 1 in absolute magnitude, then the map (17) is a contraction, and the point (18) is stable. Proof is over.

The $\Delta t$ increases the only one variant of stability loss is able. It is exit of the eigenvalue $\rho_1$ at $\Delta t = 1/\tilde{\lambda}$ beyond -1. The Hopf bifurcation takes place in this case:

$$\begin{pmatrix} Q_{i+4} \\ \lambda_{i+4} \end{pmatrix} = \begin{pmatrix} Q_i \\ \lambda_i \end{pmatrix}, \begin{pmatrix} Q_{i+6} \\ \lambda_{i+6} \end{pmatrix} = \begin{pmatrix} Q_{i+2} \\ \lambda_{i+2} \end{pmatrix}.$$

So, it is of interest: is this bifurcation subcritical or supercritical and is transfer to a strange attractor possible. Though, the following proposition gives an important information.

Proposition. There exists at least one attracting set of the map (17) at any $\Delta t > 1/\tilde{\lambda}$. All the attracting sets are inside the rectangle

$$Q \in [Q_{min}, Q_{max}], Q_{max} = \frac{\tilde{Q}}{a}\exp(a-1), Q_{min} = Q_{max}\exp(a(1-Q_{max}/\tilde{Q})), a = 2\Delta t\tilde{\lambda}, \quad (20)$$

$$\lambda \in [\lambda_{min}, \lambda_{max}], \lambda_{min} = \tilde{\lambda}(1 - Q_{max}/\tilde{Q}) + \delta_t\Delta t, \lambda_{max} = \tilde{\lambda}(1 - Q_{min}/\tilde{Q}) + \delta_t\Delta t.$$

Any initial value for the process (17) lays in a region of attraction of one of these sets.

Proof. Taking into account that $F_Q$, nor $F_\lambda$ do not depend on $\lambda$ it is enough to study the map

$$q_{i+2} = q_i \exp(a(1-q_i)) \quad (21)$$

at $a > 2$, where $q_i = Q_i/\tilde{Q}$. A fixed point of the map is $\bar{q} = 1$. The function of the map has a maximum in the point $1/a$:

$$q_{max} = \frac{\exp(a-1)}{a}. \quad (22)$$

If $q_i$ equals $q_{max}$, then

$$q_{i+2} = q_{min} = \frac{\exp(a-1)}{a}\exp(a - \exp(a-1)). \quad (23)$$

If $q_i \in [q_{min}, 1]$, then $q_i \leq q_{i+2} \leq q_{max}$. If $q_i \in [1, q_{max}]$, then $q_{min} \leq q_{i+2} \leq 1$. So, if $q_i \in [q_{min}, q_{max}]$, then $q_{i+2} \in [q_{min}, q_{max}]$ too. If $q_i > q_{max}$, then $q_{i+2} < q_{min}$. If $q_i < q_{min}$, then the further values of $q$ grow till some next $q_i$ is in the segment $[q_{min}, q_{max}]$. Reversing to the variables $Q$ and $\lambda$, we obtain the proposition rectangle. Proof is over.

Proposition. If the stability condition (19) is true, the process (21) converges to the fixed point (18) from any initial point.

Proof. The proof is similar to the previous one, but after $q_i \in [q_{min}, 1]$ the iterations converge to the fixed point. Proof is over.

The dependences (22) and (23) are given at Fig. 1 to represent a pictorial frame of maximal deviation from the fixed point.

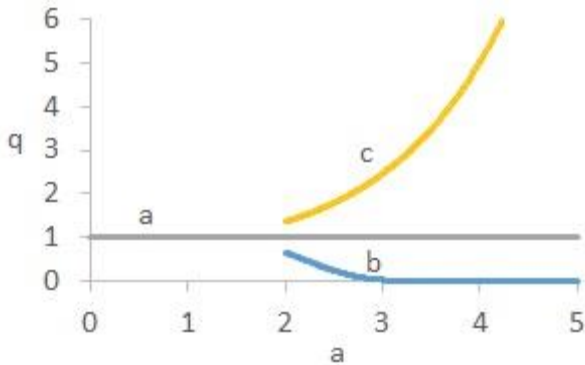

Fig. 1. The fixed point $\bar{q}$ (line a) and the functions $q_{min}(a)$ (curve b) and $q_{max}(a)$ (curve c).

Fig. 2a represents converging of the process to the fixed point at $a < 2$, i.e. before the stability condition (19) is violated. The Hopf bifurcation at $a = 2$ is supercritical. Thus, the amplitude of oscillations is much less firstly than the restriction (20) predicts. Converging to the periodic solution after this bifurcation is on the Fig. 2b. The period doubling bifurcation occurs at $a \approx 2.5265$. Converging to the cycle with doubled period is on the Fig. 2c.

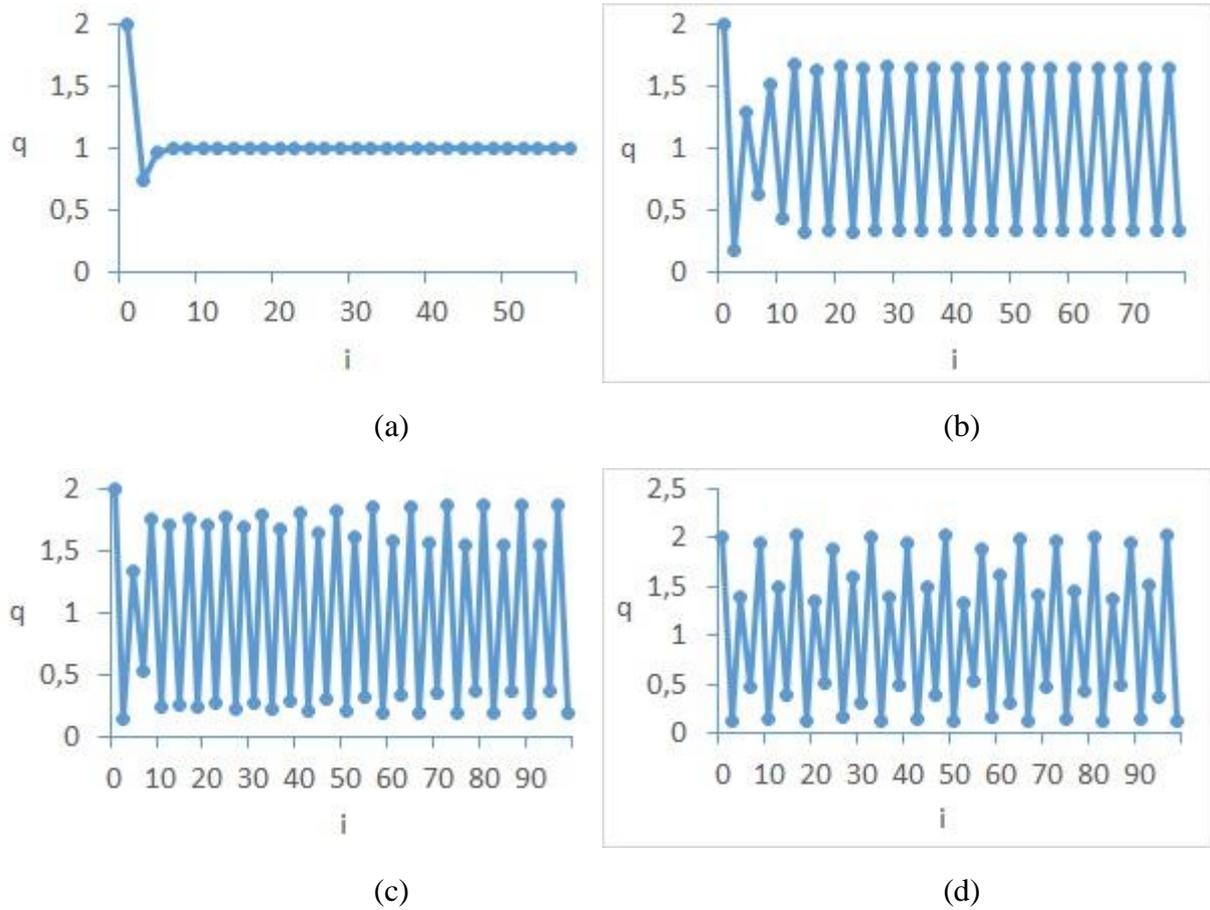

(a)  (b)  (c)  (d)

Fig. 2. a) Tendency of the process (21) to the steady state at $a = 1$. b) Tendency of the process (21) to the cycle at $a = 2.4$. c) Tendency of the process (21) to the cycle at $a = 2.58$. d) Regime of a strange attractor at $a = 2.7$. $L \approx 0.09$.

Then a sequence of period doublings takes place: $a \approx 2.5265, 2.6564, 2.6846, 2.6907, \ldots$ Approximately at $a \approx 2.6924$ a strange attractor appears. An example of the algorithm functioning in a strange attractor regime is on the Fig. 2d. The positive Lyapunov exponent justifies that the process is really chaotic. It is easy to see from Fig. 1 and Fig. 2d, that the restriction (20) predicts maximal and minimal values of $q$ in this case well. Different effects can be observed at different values of $a$ and a larger number of the process iterations: beats, temporary (deceptive) exit to steady state etc. Though, all this takes place inside the rectangle (20).

The map (17) is such that $F_Q$ and $F_\lambda$ do not depend on $\lambda$. Thus, all the points $(Q_{i+2}, \lambda_{i+2})$ lay on one and the same curve:

$$Q = \tilde{Q}\left(1 - \frac{\lambda - \delta_t \Delta t}{\tilde{\lambda}}\right) \exp(2\Delta t(\lambda - \delta_t \Delta t)). \tag{24}$$

For example, a strange attractor, corresponding to Fig. 2d, is pictured at Fig. 3 in the plane $(q, \lambda)$.

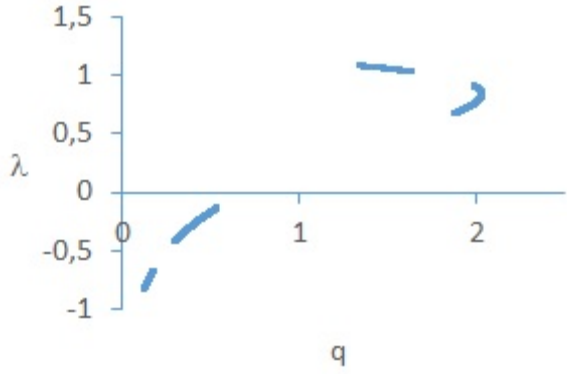

Fig. 3. $10^5$ iterations of the algorithm 1. $\tilde{\lambda}=1$, $\delta_t=0.1534$, $\Delta t=1.35$. Looking as solid parts of the curve (24) are sets of separate point really.

## 5. Bifurcation properties of the algorithm 2 on the equation (6).

The following formulas determine the map (14) in the case of the algorithm 2

$$Q_{i+1} = Q_i \frac{\exp(\delta_N(N_i+\Delta N)(t_i+\Delta t)+\delta_t(t_i+\Delta t)^2/2)}{\exp(\delta_N(N_i+\Delta N)t_i+\delta_t t_i^2/2)} =$$
$$= Q_i \exp(\lambda_i \Delta t + \delta_N \Delta N \Delta t + \delta_t \Delta t^2/2), \qquad (25)$$
$$\lambda_{i+1} = \lambda_i + \delta_N \Delta N + \delta_t \Delta t,$$

where $\Delta N$ can be found by formulas (9), (13). The map (25) has a unique fixed point

$$\overline{Q} = \tilde{Q}\left(1 - \frac{\delta_t \Delta t(1/2 - 1/a)}{\tilde{\lambda}}\right),$$
$$\overline{\lambda} = \delta_t \Delta t / 2, \qquad (26)$$

where $a = \delta_N / \tilde{\delta}_N > 0$. The formulas (26) are conscious if $\overline{Q}>0$ and $\tilde{Q}>0$.

<u>Proposition.</u> It is necessary and sufficient for stability of the steady state $\{\overline{Q},\overline{\lambda}\}$ for all small enough $\Delta t>0$, the following conditions to be true

$$\begin{cases} a<2, \\ ab>1, \end{cases} \qquad (27)$$

where $b = \tilde{\lambda}/\tilde{Q}$.

<u>Proof.</u> Since the formulas (26) a conscious at $\overline{Q}\tilde{Q}>0$, the following condition should be fulfilled

$$\tilde{\lambda} > \delta_t \Delta t(1/2 - 1/a).$$

It is fulfilled in two cases:

$$\begin{bmatrix} a) & \delta_t > 0, & a > 2 \\ b) & \delta_t < 0, & 0 < a < 2 \end{bmatrix} : \quad \Delta t \in (0, \frac{\tilde{\lambda}}{\delta_t (1/2 - 1/a)}), \qquad (28)$$

$$c) \quad \textit{other cases}: \quad \Delta t \in (0, \infty).$$

In particular, it is always fulfilled for small enough $\Delta t$.

The eigenvalues of the Jacoby matrix of the map (25) in the point (26) are

$$\theta_{1,2} = \frac{(2-a) - h \cdot ab \pm \sqrt{a(h^2 \cdot ab^2 + h \cdot 2b(a-2) + a)}}{2}, \qquad (29)$$

where $h = \overline{Q} \Delta t$. Only $h = h(\Delta t)$ depends on $\Delta t$ in (29), moreover, $h(0) = 0$ and $h'(0) > 0$. If $h = 0$, the eigenvalues of (29) equal $\theta_1 = 1$ and $\theta_2 = 1 - a$. The second condition in (27) means negative derivative of $\theta_1$ on $h$ at $h = 0$ (an analogous variant $a = 2$ and positive derivative of $\theta_2$ is impossible). Proof is over.

<u>Proposition.</u> Let the radical expression in the formula

$$\Delta t_0 = \frac{b\left(\tilde{Q} - \sqrt{\tilde{Q}^2 - 8\delta_t (0.5 - 1/a)(2-a)/(ab^2)}\right)}{2\delta_t (0.5 - 1/a)}. \qquad (30)$$

is positive. Then the point (26) loses stability for increasing $\Delta t$ then and only then when $\Delta t$ overcomes the value $\Delta t_0$. The eigenvalue $\theta_2$ passes over the value -1 in this case. If the radical expression is negative, then the point (26) is stable for any $\Delta t$ of the condition (28) b).

<u>Proof.</u> At $a \geq 1$ the radical expression in (29) is nonnegative and both eigenvalues are real. At $a \in (0,1)$ there is an interval on the positive semiaxis

$$h \in (h_1, h_2), \quad h_1 = \frac{2 - a - 2\sqrt{1-a}}{ab} > 0, \quad h_2 = \frac{2 - a + 2\sqrt{1-a}}{ab} > h_1,$$

in which the roots $\theta_1$ and $\theta_2$ are complex conjugate.

If the discriminant in (29) is greater than zero, then, subject to the restrictions for $a$, $b$ and $h$, the eigenvalue $\theta_1$ can't equal 1, and the eigenvalue $\theta_2$ equals -1 then and only then, when $h = h_0$:

$$h_0 = \frac{2(2-a)}{ab}.$$

If the radical expression in (29) is less than zero, then, subject to the restrictions for $a$, $b$ and $h$, the eigenvalues $\theta_1$ and $\theta_2$ can't be equal to 1 in absolute magnitude. Note, that $h_0 > h_2$. Considering $h = h(\Delta t)$ and relations (27), (28) we obtain the condition (30). Proof is over.

It is clear from the propositions about stability and its loss that the algorithm 2 reveals much more complex properties. The steady sates depend on greater number of parameters, moreover, $\overline{Q} \neq \tilde{Q}$; there are two, but not one, conditions of stability; and more variants of losing or keeping stability for $\Delta t$ increasing. Besides, in the case of algorithm 2, if the conditions

$$\begin{cases} \delta_t < 0, \\ a < 2\delta_t \Delta t /(\delta_t \Delta t - 2\tilde{\lambda}), \end{cases} \quad (31)$$

are fulfilled, then the steady state

$$\begin{cases} \overline{Q} = 0, \\ \overline{\lambda} = \tilde{\lambda} + \delta_t \Delta t / a. \end{cases} \quad (32)$$

is stable. This fixed point can play essential role even in the case if it is unstable. The exponent is present in the formula for $Q_{i+1}$ in (25), thus, the value $Q_{i+1}$ can prove quite close to 0 and, so, can remain there for quite a long time. For a numerical investigation, it can turn to computer zero. After this, $\{Q_{i+1}\}$ is changeless, and $\{\lambda_{i+1}\}$ goes to $\overline{\lambda}$ (see. (32)).

The following values of parameters were chosen for numerical investigation of the algorithm 2 on the equation (6): $\delta_N=0.2$, $\delta_t=0.25$, $N_0=1$, $\tilde{Q}=1$, $\tilde{\lambda}=4$, $\tilde{\delta}_N=0.5$. The Fig. 4a represents tending to the stable steady state when the conditions (27) are true.

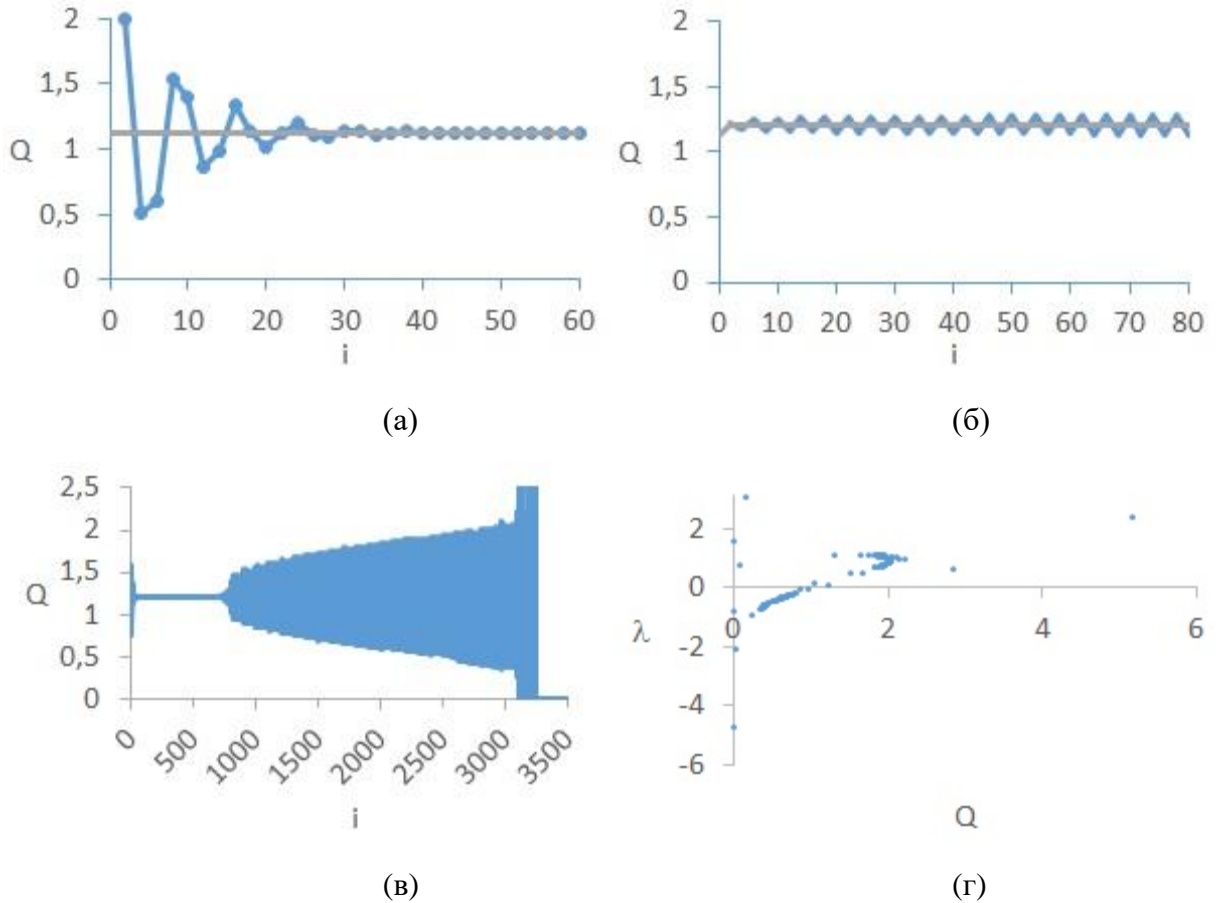

(а)　　　　　　　　　　　　　　　(б)

(в)　　　　　　　　　　　　　　　(г)

Fig. 4. a) Tending to the stable fixed point (26) at $\Delta t = 1$. $Q_0=2$, $\lambda_0 = \delta_N N_0$. b) Tending to the periodic solution at $\Delta t = 1.66$. $Q_0 = 1.01 \cdot \overline{Q}$, $\lambda_0 = 0.99 \cdot \overline{\lambda}$. c) A sequence of $Q_i$ values. The initial value $\Delta t = 1.6$, then, $\Delta t$ enlarges with 0.01 every 100 steps. d) A strange attractor, hypothetically, at $\Delta t = 1.89$.

The Hopf bifurcation takes place at $\Delta t = \Delta t_0 \approx 1.65685$, and the fixed point becomes unstable. Converging to the appeared periodical solution after this bifurcation is drawn on the Fig. 4b.

The initial condition is near the steady state since the region of convergence is uncertain. Period doublings takes place when $\Delta t$ continues increasing: 1.84664, 1.87434, 1.87924, 1.88027, 1.88049, ... Strange attractor appears approximately at $\Delta t \approx 1.8805$. The sequence of period doublings and transfer to the strange attractor, hypothetically, is on the Fig. 4c. It is said "hypothetically" here because calculation of the Lyapunov exponents is strongly obstructed by turning of $Q_i$ to computer zero. Some points are out the picture scale. The are segments where a long sequences of $Q_i$ is very close to zero, then $Q_i \approx 0.5$, then $Q_i \approx 1470$ (this point is out of the scale) and then $Q_i = 0$ for computer precision. A sequence of approximately 135 points on the Fig. 4d hypothetically corresponds to a strange attractor. Several points are out of the picture scale. The maximal Lyapunov exponent is uncalculated because $Q_i$ turned to computer zero.

### 6. Discussion of results.

Comparison of the algorithms 1 and 2 on the equation (6) shows substantial advantage of the algorithm 1. The steady state of the value $Q$ of this algorithm is $\tilde{Q}$ (18) if the map (15) begins with an even step. Moreover, it is important that this value does not depend on the parameters and the system state, but depends only on the parameter of the algorithm itself. The same relates to the stability condition (19) of this steady state, i.e. the only parameters of the algorithm determine it. The steady state of the value $Q$ (26) can strongly differ from $\tilde{Q}$ in the case of the algorithm 2. Moreover, it depends on $\delta_t$ besides the parameters of the algorithm. The $\delta_N$ does not influences the steady states and stability conditions for both algorithms. There is convergence to the fixed point (18) from any initial value, if this fixed point is stable, in the case of the algorithm 1. The stability condition is always true for the algorithm 1, if $\Delta t$ is small enough. The stability condition (27) of the algorithm 2 is always true for small enough $\Delta t$, if the parameters $\tilde{\delta}_N$ and $\tilde{\lambda}$ are proper. It is always possible to choose these parameters properly. One can choose a large enough $\tilde{\delta}_N$ to satisfy the condition $a < 2$, and then a large enough $\tilde{\lambda}$ to satisfy the condition $ab > 1$.

The step value $\Delta t$ growth, every algorithm leads to Hopf bifurcation and the fixed point stability loss. The attracting sets are always inside the rectangle (20) for the algorithm 1. Thus, even the work in the strange attractor regime is predictable enough, although the parameters of the algorithm are still worth choosing for the steady state (18) be stable. Forecasting for the algorithm 2 is much less, what is particularly because of the special role of the steady state (32). The sequences $\{(Q_i, \lambda_i)\}$ can get in the vicinity of this steady state and be there for a long time even if this steady state is unstable.

The said above relates to control of the equation (6) by idealized algorithms. The distinctions of the real system from the equation (6) are possible. The difference measuring of the value $\lambda$ and nonzero interval of the controlling value $N$ changing are possible too. All this distinctions from the situation investigated tend to zero for tending to zero of the value $\Delta t$. So, for small enough $\Delta t$, the steady states of the real system will tend to the steady states (18) and (26), and the stability conditions will tend to the stability conditions (19) and (27). Moreover, if the system properties differ only a little from investigated case, then the results will be similar for the similar value $\Delta t$.

I.e., the results will be similar numerically also. Taking into account possibility of small $\Delta t$ choosing the algorithms seem to be applicable to many systems: to wind generators, for example [11-12].

The two algorithms for controlling solutions of the systems of differential equations are investigated in the article. An equation is constructed under assumption that the systems keep or almost keep some definite properties, and this equation can be the main equation for investigation of the idealized algorithms properties. It is so because the method of investigation implies some suppositions fulfilled for this equation. It is possible to investigate the differences appearing in the real systems separately, and it is out of this article scope. The differences tend to zero for tending to zero the step of control. The fixed points and stability conditions are found for the both algorithms. The bifurcations for violation of the stability conditions are investigated. One of the algorithms acts much more reliably (in the case of the violation of the stability condition as well). It may be recommended for practical use. Another algorithm is of much more complex properties. Exit beyond the stability conditions leads to periodic solutions and strange attractors for both algorithms.